\documentclass{amsart}

\usepackage{amssymb, amsmath}

\newcommand{\N}{\mathbb{N}}  
\newcommand{\Z}{\mathbb{Z}}  
\newcommand{\ZZ}[1]{{\Z/{#1}\Z}}  
\newcommand{\R}{\mathbb{R}}  
\newcommand{\CC}{\mathbb{C}}  
\newcommand{\T}{\mathbb{T}}  

\newcommand{\fr}{\mathfrak}  

\newcommand{\goth}{\mathfrak}
\newcommand{\mc}{\mathcal}  

\newtheorem{Th}{Theorem}

\newtheorem{Lem}[Th]{Lemma}

\begin{document}
\numberwithin{Th}{section} 
\setcounter{section}{-1} 

\numberwithin{equation}{section} 

\title{Ergodic Abelian actions with homogeneous spectrum}

\author{Alexandre I. Danilenko}
\address{Institute for Low Temperature Physics \& Engineering of National Academy of Sciences of Ukraine, 47 Lenin Ave., Kharkov, 61164, UKRAINE}
\email{alexandre.danilenko@gmail.com}

\author{Anton V. Solomko}
\email{solomko.anton@gmail.com}
\thanks{The second named author was partially supported by Akhiezer fund}
\maketitle

\begin{abstract}
It is shown that for each $N>0$ and for a wide class of Abelian non-compact locally compact  second countable groups $G$ including
all infinite countable discrete ones and $\Bbb R^{d_1}\times\Bbb Z^{d_2}$ with $d_1,d_2\ge 0$, there exists a weakly mixing
probability preserving  $G$-action with a homogeneous spectrum of
multiplicity $N$.
\end{abstract}

\section{Introduction}

Assume that $G$ is a non-compact locally compact  second countable Abelian group.
Let $T=(T_g)_{g\in G}$ be a measure preserving action of $G$ on a standard probability space $(X,\goth B,\mu)$.
The corresponding unitary representation $U_T$ of $G$ in $L^2(X,\mu)$,
$$
U_T(g)f:=f\circ T_g^{-1},
$$
is called the Koopman representation of $G$.
As usual, we will  consider $U_T$ only on $L^2_0(X,\mu):=L^2(X,\mu)\ominus\Bbb C$,
that is on the subspace of zero mean functions.
It is easy to construct ergodic $G$-actions $T$ such that $U_T$ has a  simple spectrum.
For instance, any action with a purely discrete spectrum or, more generally, any rank-one action has a simple spectrum.
In the present paper we consider
\begin{enumerate}
\item[({\bf Pr})]
 {\it the problem of existence  of ergodic $G$-actions
whose Koopman representations have a homogeneous spectrum of multiplicity greater than $1$}.
\end{enumerate}

In the case $G=\Bbb Z$ this problem is attributed to Rokhlin.
The first examples of ergodic transformations with a homogeneous spectrum of multiplicity 2 appeared in 1999 in \cite{Ry1} and \cite{A1}.
In 2005, Ageev \cite{A2} proved the existence of ergodic transformations with a homogeneous spectrum of arbitrary multiplicity (see also subsequent works  \cite{D1} and \cite{Ry2} for simplifications and concrete examples of such transformations).
In the case $G=\Bbb Z^2$, Konev explains in his recent work \cite{Ko} how to achieve  homogeneous spectrum of multiplicity 2 for  ergodic $G$-actions.
We also note that Ryzhikov stated~(Pr) for $G=\Bbb R$  in 2006 (see problem (P9) in \cite[\S5]{D2}).

We now state the main result of the paper.
Recall that each Abelian locally compact Polish group $G$ is isomorphic to the product $\Bbb R^p\times G'$, where $G'$ is a locally compact group containing an open compact subgroup $G'_0$ \cite[Theorem~24.30]{HR}.

\begin{Th}\label{MainTh}
Let one of the following be satisfied:
\begin{enumerate}
\item[(i)]
$p>0$,
\item[(ii)]
$p=0$ and $G_0'$ is a direct summand in $G_0$,
\item[(iii)]
$p=0$, $G_0'$ is not a direct summand in $G_0$ and there is no $k>0$ such that $k\cdot g=0$ for all $g\in G'/G_0'$.
\end{enumerate}
Then for each $N>1$, there exists a  weakly mixing  $G$-action  with a homogeneous spectrum of multiplicity $N$.
\end{Th}

In particular, the conditions of the theorem hold if $G$ is any infinite countable discrete Abelian group or $G=\Bbb R^p$ or  the direct product of a discrete group with $\Bbb R^p$ or if $G$ is a field of $q$-adic numbers, etc.

To show this theorem we  use   a ``generic'' argument as in \cite{A2}, \cite{D1} and \cite{Ry2}:  given $G$ and $N$, we introduce an auxiliary
non-Abelian group, topologize the set of all measure preserving actions of this group and show that certain two properties of these actions are generic. Hence there is an action  possessing both of them. We then show that the two properties imply the homogeneous spectrum of  multiplicity  $N$.
Even in the case of $G=\Bbb Z$, our auxiliary group  is different from those considered in \cite{A2}, \cite{D1} and \cite{Ry2}. It  includes $G$ as a direct summand.

We thank V.~Ryzhikov for his stimulating questions.

\section{Countable discrete Abelian group actions with homogeneous spectrum} \label{DiscSec}

In this section we prove Theorem~\ref{MainTh} in the case when $G$ is discrete, countable, infinite and Abelian. In general, the proof goes along the lines developed in \cite{D1}. However, the crucial step is to choose a new auxiliary group $\Gamma$ in a ``right'' way (see below).

Denote by $\mc U(\mc H)$ the group of unitary operators on a separable Hilbert space $\mc H$.
We endow $\mc U(\mc H)$ with the (Polish) strong operator topology (which on $\mc U(\mc H)$ is also the weak operator topology).
Given a discrete countable group $F$, we furnish the product space $\mc U(\mc H)^F$ with the product topology and denote by $\mc U_F(\mc H)\subset\mc U(\mc H)^F$ the subset of all unitary representations of $F$ in $\mc H$.
Obviously, $\mc U_F(\mc H)$ is closed in $\mc U(\mc H)^F$ and hence Polish in the induced topology.

Given a standard non-atomic probability space $(X,\fr B,\mu)$, let $\text{Aut}(X,\mu)$ stand for the group of invertible $\mu$-preserving transformations of $X$.
By an \emph{action} $T$ of $F$ we mean a group homomorphism $T\colon F\ni f\mapsto T_f\in\text{Aut}(X,\mu)$.

Denote by $\mc A_F\subset\text{Aut}(X,\mu)^F$ the subset of all measure-preserving actions of $F$ on $(X,\fr B,\mu)$.
Recall that $U_T$ denotes the Koopman representation of $F$ associated with $T\in\mc A_F$.
We endow $\mc A_F$ with the weakest topology which makes continuous the mapping
$$
\mc A_F\ni T \mapsto U_T \in \mc U_F(L_0^2(X,\mu)).
$$
It is Polish.
 It is easy to verify that a sequence $T^{(m)}$ of $F$-actions converges to $T$ if and only if $\mu(T^{(m)}_f A\bigtriangleup T_f A)\to 0$ as $m\to\infty$ for each $f\in F$  and $A\in \fr B$.
There exists a natural continuous action of $\text{Aut}(X,\mu)$ on $\mc A_F$ by conjugation:
$$
(R\cdot T)_f = R T_f R^{-1} \quad \text{for $R\in\text{Aut}(X,\mu)$, $T\in\mc A_F$, $f\in F$.}
$$

The following three lemmata are well known in the case of $\Bbb Z$-actions. Since we were unable to find in the literature their generalization to the general Abelian case, we provide these lemmata  with short proofs.

\begin{Lem}\label{WM}
The subset of weakly mixing $G$-actions is $G_\delta$ in $\mc A_G$.
\end{Lem}

\begin{proof}
First let us show that the subset $\mathcal E := \{ T\in \mc A_G \mid T~\text{is ergodic} \}$ is $G_\delta$ in $\mc A_G$.
Given $T\in\mc A_G$, denote by $P_\text{inv}$ the orthogonal  projection from $L_0^2(X,\mu)$ to the subspace of  $T$-invariant functions in $L_0^2(X,\mu)$.
It is well known that $T$ is ergodic if and only if $P_\text{inv}=0$.
Fix a F{\o}lner sequence $\{F_n\}_{n\in\N}$ in $G$.
It follows from the von Neumann theorem \cite[Theorem~3.33]{Gl} that for each  $T\in\mc A_G$,
$$
\lim_{n\to\infty}\frac{1}{|F_n|}\sum_{g\in F_n}U_T(g)\varphi=P_\text{inv}\varphi
$$
for every $\varphi\in L_0^2(X,\mu)$.
Therefore, $T$ is ergodic if and only if for any $\varepsilon>0$, $\varphi\in L_0^2(X,\mu)$ and $N\in\N$, there exists $n>N$ such that $\bigl\| \frac{1}{|F_n|}\sum_{g\in F_n}U_T(g) \varphi \bigr\| < \varepsilon$.
Let $\{\varphi_k\}_{k\in\N}$ be a dense subset in $L_0^2(X,\mu)$.
Then $\mc E$ is $G_\delta$  in $\mc A_G$ because
$$
\mathcal E = \bigcap_{m=1}^\infty \bigcap_{k=1}^\infty \bigcap_{N=1}^\infty \bigcup_{n=N}^\infty
\Bigl\{ T\in \mc A_G \mid \Bigl\| \frac{1}{|F_n|}\sum_{g\in F_n}U_T(g) \varphi_k \Bigr\| < \frac{1}{m} \Bigr\}.
$$
This and the following facts:
\renewcommand{\labelenumi}{(\theenumi)}
\renewcommand{\theenumi}{\roman{enumi}}
\begin{enumerate}
\item $T$ is weakly mixing if and only if $T\times T$ is ergodic;
\item the map $T\mapsto T\times T$ is continuous,
\end{enumerate}
imply that the subset of weakly mixing $G$-actions is also $G_\delta$ in $\mc A_G$.
\end{proof}

Let $U$ be a unitary representation of $G$ in a separable Hilbert space $\mc H$.
By the cyclic subspace of $\psi\in\mc H$ we mean the smallest closed $U$-invariant subspace $\mc H(\psi)$ of $\mc H$ containing $\psi$.
Let also $d(\varphi,\mc H')$ stand for the distance from a vector $\varphi\in \mc H$ to a subspace $\mc H'$ in $\mc H$.

By the spectral theorem for $U$, there exist a probability measure $\sigma$ on the compact Polish dual group $\widehat{G}$ and a Borel  field of Hilbert spaces $\widehat{G}\ni\omega\mapsto\mc H_\omega$
such that (up to unitary equivalence)
\begin{equation}\label{Hdecomp1}
\mc H = \int_{\widehat{G}}^\oplus \mc H_\omega d\sigma(\omega)
\quad \text{and}\quad
U_g = \int_{\widehat{G}}^\oplus \omega(g)I_\omega d\sigma(\omega)
\text{ for each }g\in G,
\end{equation}
where 
$I_\omega$ is the identity operator on $\mc H_\omega$.
Moreover, for every $\psi\in\mc H$
\begin{equation}\label{CycSpDecomp}
\mc H(\psi) = \int_{\widehat{G}}^\oplus \mc H'_\omega d\sigma(\omega),
\end{equation}
where $H'_\omega$ is a subspace in $H_\omega$, $\dim \mc H'_\omega=0$ or 1 for $\sigma$-almost all $\omega\in\widehat{G}$.
The map $k\colon\widehat{G}\to\N\cup\{\infty\}$, $k(\omega) = \dim \mc H_\omega$ is called the \emph{multiplicity function} of $U$.
Denote by $\mc M(U)$ the set of essential values for the multiplicity function of $U$.
By definition, $U$ has \emph{a homogeneous  spectrum of multiplicity $N$} if $\mc M(U)=\{N\}$, and $U$ has \emph{ a simple spectrum} if $\mc M(U)=\{1\}$.
Given $T\in\mc A_G$, let $\mc M(T):=\mc M(U_T)$.

\begin{Lem}[Katok-Stepin lemma for Abelian actions]\label{KS}
If $U\in\mc U_G(\mc H)$ does not have simple spectrum then there exist unit vectors $\varphi_1,\varphi_2\in\mc H$ such that for all $\psi\in\mc H$
$$
d^2(\varphi_1,\mc H(\psi))+d^2(\varphi_2,\mc H(\psi))\geqslant 1.
$$
\end{Lem}

\begin{proof}
By the spectral theorem for $U$, the decomposition (\ref{Hdecomp1}) holds.
By the condition of the lemma, there is a subset $A$  in $\widehat{G}$, $\sigma(A)>0$, such that $k(\omega)\geqslant 2$ for each $\omega\in A$.
 Then we can choose measurable fields $\widehat{G}\ni\omega\mapsto \varphi_1(\omega), \varphi_2(\omega)\in \mc H_\omega$ in such a way that $\|\varphi_1(\omega)\|=\|\varphi_2(\omega)\|=\frac{1}{\sqrt{\sigma(A)}}$, $\varphi_1(\omega)\perp \varphi_2(\omega)$ for $\omega\in A$, and $\varphi_1(\omega)=\varphi_2(\omega)=0$ for $\omega\in\widehat{G}\setminus A$.
Let $\psi\in\mc H$.
It follows from (\ref{CycSpDecomp}) that
$$
\mc H(\psi) = \int_{\widehat{G}}^\oplus \mc H'_\omega d\sigma(\omega) \quad\text{with}\quad \dim \mc H'_\omega \leqslant 1.
$$
Applying Bessel inequality we obtain $d^2(\varphi_1(\omega),\mc H'_\omega)+d^2(\varphi_2(\omega),\mc H'_\omega)\geqslant \frac{1}{\sqrt{\sigma(A)}}$ for any $\omega\in A$.
Hence
\begin{multline*}
d^2(\varphi_1,\mc H(\psi))+d^2(\varphi_2,\mc H(\psi))=
\int_A \bigl( d^2(\varphi_1(\omega),\mc H'_\omega)+d^2(\varphi_2(\omega),\mc H'_\omega) \bigr) d\sigma(\omega) \geqslant 1.
\end{multline*}
\end{proof}

\begin{Lem}\label{SimSp}
The subset of actions with a simple spectrum is $G_\delta$ in $\mc A_G$.
\end{Lem}

\begin{proof}
By the definition of topology on $\mc A_G$, it suffices to show that
 given a separable Hilbert space $\mc H$, the subset $\mc S\subset\mc U_G(\mc H)$ of unitary representations with a simple spectrum is $G_\delta$.
Fix a dense subset $\mc D=\{\varphi_i\}_{i\in\N}$  in $\mc H$.
Enumerate the elements of $G$ as $\{g_j\}_{j\in\N}$.
Define a subset $\mc S_1\subset\mc U_G(\mc H)$ as follows: $U\in\mc S_1$ if and only if given $\varphi_1,\ldots,\varphi_n\in\mc D$ and $\varepsilon>0$, there exist $\psi\in\mc H$, $N\in\N$ and $\lambda_{i,1},\ldots,\lambda_{i,N}\in\CC$, $i=1,\ldots,n$, such that
$$
\Bigl\| \varphi_i - \sum_{j=1}^N \lambda_{i,j} U(g_j)\psi \Bigr\| < \varepsilon, \quad i=1,\ldots,n.
$$
We see that $\mc S_1$ is $G_\delta$  in $\mc U_G(\mc H)$:
$$
\mc S_1 = \bigcap_{n=1}^\infty \bigcap_{m=1}^\infty \bigcup_{N=1}^\infty \bigcup_{\psi\in\mc H} \bigcup_{\{\lambda_{i,j}\}\subset\CC}\{ U\in\mc U_G \mid \max_{1\leqslant i\leqslant n} \Bigl\| \varphi_i - \sum_{j=1}^N \lambda_{i,j} U(g_j)\psi \Bigr\| < \frac{1}{m} \}.
$$
If $U$ has a simple spectrum then there is $\psi\in\mc H$ with $\mc H(\psi)=\mc H$. Therefore, $\mc S_1\supset\mc S$.
The converse inclusion  follows from Lemma~\ref{KS}:
if $U\notin\mc S$ then we can choose unit vectors $\varphi_1,\varphi_2\in\mc H$ in such a way that $d^2(\varphi_1,\mc H(\psi))+d^2(\varphi_2,\mc H(\psi))\geqslant 1$ for each $\psi\in\mc H$, and hence $U\notin\mc S_1$.
\end{proof}

Suppose we are given a countable discrete Abelian group $J$.
Then we define a `cyclic' group automorphism $A\colon J^\ZZ{N}\to J^\ZZ{N}$ by setting
$$
(Aj)(i):=j(i-1), \quad j\in J^\ZZ{N}, \quad i\in\ZZ{N},
$$
and denote by $\Gamma$ a semidirect product
$
G\times J^\ZZ{N}\rtimes_A(\ZZ{N})
$
with the multiplication law as follows:
$$
(g,j,n)(h,k,m):=(g+h,j+A^n k,n+m),\quad g,h\in G,\: j,k\in J^\ZZ{N},\: n,m\in \ZZ{N}.
$$

We now isolate an Abelian subgroup
$
\Lambda :=\{ (g,0,n) \mid g\in G, n\in \ZZ{N} \}
$
in $\Gamma$.

The following lemma is crucial in the proof of the main theorem.
We show how to choose $J$ to construct a free $\Gamma$-action whose restriction to $\Lambda$ has a simple spectrum.

\begin{Lem} \label{MainLem}
There exists $J$ and a free action $T$ of the corresponding group $\Gamma$ such that $T\upharpoonright\Lambda$ is ergodic with a pure point spectrum.
\end{Lem}

\begin{proof}
Let $K$ be a compact second countable Abelian group and let $R\colon G\to K$ be a  one-to-one homomorphism with $\overline{R(G)}=K$.
To see that such $K$ and $R$ exist, let $\{\chi_n\}_{n\in\N}\subset\widehat{G}$ be such that for each $0\neq g\in G$ there is $n$ with $\chi_n(g)\neq 1$.
Define a continuous embedding $R\colon G\to \T^\N$ by setting $R(g):=(\chi_1(g),\chi_2(g),\ldots)$.
It remains to set $K:=\overline{R(G)}$.

We claim that there are a countable discrete Abelian group $J$ and one-to-one homomorphisms $S_i\colon J\to K$, $i\in\ZZ{N}$, such that
the subgroups
\begin{equation}\label{Idepend}
R(G), S_0(J),\dots,S_{N-1}(J) \text{ are independent in } K,
\end{equation}
i.~e., $R(g)+S_0(j_0)+\cdots+S_{N-1}(j_{N-1})=0$ implies $g=0$, $j_0=\cdots=j_{N-1}=0$.
For that we consider separately two cases.
Denote by $K_0$ the torsion subgroup of $K$.

\emph{First case}.
Suppose that $K_0$ is meager.
Given a countable subgroup $L\subset K$, let $\widetilde{L}:=\{k\in K \mid \exists m>0: mk\in L \}$.
Obviously, $\widetilde{L}\supset L\cup K_0$.
Since the quotient group $\widetilde{L}/K_0$ is countable and $K_0$ is meager, $\widetilde{L}$ is also meager.
The Baire category theorem implies that $K\setminus \widetilde{L}\neq\emptyset$, i.~e., there exists an element $0\neq k\in K$ of infinite order such that $\langle k\rangle\cap L=\{0\}$.
We will use this fact repeatedly.
First put $L:=R(G)$ and choose an element $k_0\in K\setminus K_0$ with $\langle k_0\rangle\cap R(G)=\{0\}$.
Second put $L:=R(G)+\langle k_0\rangle$ and choose an element $k_1\in K\setminus K_0$ with $\langle k_1\rangle\cap (R(G)+\langle k_0\rangle)=\{0\}$.
Continuing in the same way, we construct elements $k_0,\ldots,k_{N-1}\in K$ of infinite order such that the subgroups $R(G),\langle k_0\rangle,\ldots\langle k_{N-1}\rangle$ are independent.
 It remains to set $J:=\Z$ and define embeddings $S_i\colon J\to K$, $i\in\ZZ{N}$, by  $S_i(j):=jk_i$.
Clearly, (\ref{Idepend}) is satisfied for  $S_i$,  $i\in\ZZ{N}$.

\emph{Second case}.
Now suppose that $K_0$ is nonmeager.
Then by \cite[Theorem~2.3]{VTCh}, $K_0-K_0$ contains a neighborhood of the identity.
Thus, $K_0$ is open in $K$.
This and the compactness of $K$ implies $K=K_0$.
It follows from \cite[Theorem 25.9]{HR} that $K$ splits into direct product $K=(\Bbb Z/q\Bbb Z)^\N\times C$ for some  $q>1$, where $C$ is another compact Abelian group.
Enumerate the elements of $R(G)$ as $(g_m,c_m)$, $g_m\in(\ZZ{q})^\N$, $c_m\in C$, $m\in\N$.
Fix nonempty pairwise disjoint subsets $A_{i,n,m}^l, B_{i,n}$ of $\N$, $i\in\ZZ{N}$, $n,m\in\N$, $l=1,\ldots,q-1$, and
select $f_{i,n}\in(\ZZ{q})^\N$, $i\in\ZZ{N}$, $n\in\N$, in such a way that
\renewcommand{\labelenumi}{(\theenumi)}
\renewcommand{\theenumi}{\roman{enumi}}
\begin{enumerate}
\item \label{fone} $l f_{i,n}(k)\neq g_m(k)$, for $k\in A_{i,n,m}^l$,
\item \label{ftwo} $f_{i,n}(k):=1$, for $k\in B_{i,n}$,
\item \label{fthree} $f_{i,n}(k):=0$, otherwise.
\end{enumerate}
To satisfy (\ref{fone}) put $f_{i,n}(k):=0$ if $g_m(k)\neq 0$, and $f_{i,n}(k):=1$ if $g_m(k)= 0$.
Each $f_{i,n}$ is an element of period $q$ by (\ref{ftwo}).
We now set $J:=\bigoplus_{n=1}^\infty (\Z/q\Z)$ and define embeddings $S_i\colon J\to K$, $i\in\ZZ{N}$, by setting
$$
S_i((l_1,l_2,\ldots)):=l_1 (f_{i,1},0) + l_2 (f_{i,2},0) + \cdots.
$$
Then (\ref{fone}) and (\ref{fthree}) imply (\ref{Idepend}) for   $S_i$, $i\in\ZZ{N}$.

Thus the claim is proved.
Now let $X:=K\times(\ZZ{N})$.
Then $X$ is a compact Abelian group.
We  equip $X$ with Haar measure.
Now we  define an action $T$ of $\Gamma$ on $X$ by setting
$$
T_{(g,j,m)}(k,n):=(k+R(g)+S_{n+m}(j(n+m)),n+m),
$$
$(g,j,m)\in\Gamma$, $(k,n)\in X$.
It follows from (\ref{Idepend}) that $T$ is free.
To show that $T\upharpoonright\Lambda$ is ergodic and has a pure point spectrum we just note that $\Lambda$ acts on $X$ by rotations
$$
T_{(g,0,n)}x=x+(R(g),n)
$$
and the subgroup $\{ (R(g),n) \mid (g,0,n)\in\Lambda \}$ is dense in $X$.
\end{proof}

The following two lemmata are analogues of Proposition~1.2 and Lemma~1.3 from \cite{D1} respectively.
To state them we need
to introduce more subgroups in $\Gamma$:
\begin{align*}
H_i &:=\{ (0,j,0) \mid j\in J^\ZZ{N}, j(i)=-j(0), j(k)=0 \text{ for } k\neq 0, i \},\quad 1\leqslant i<N,\\
H &:=\{ (g,j,0) \mid g\in G, j\in J^\ZZ{N} \},\\
F &:=\Lambda\cap H\cong G.
\end{align*}

\begin{Lem} \label{WSresidual1}
The following two subsets are dense $G_\delta$ in $\mc A_\Gamma$:
\begin{align*}
\mc W &:=\{ T\in\mc A_\Gamma \mid T\upharpoonright H_i \text{ is weakly mixing for each $1\leqslant i<N$} \}
\text{ and}\\
\mc S &:= \{ T\in\mc A_\Gamma \mid T\upharpoonright\Lambda \text{ has a simple spectrum} \}
\end{align*}
\end{Lem}

\begin{proof}
Since for each subgroup $\Gamma_1\subset\Gamma$, the map $\mc A_\Gamma\ni T\mapsto T\upharpoonright \Gamma_1 \in\mc A_{\Gamma_1}$ is continuous,
it follows from Lemmata~\ref{WM} and \ref{SimSp} that $\mc W$ and $\mc S$ are both $G_\delta$ in $\mc A_\Gamma$.
We note that $\mc W$ and $\mc S$ are $\text{Aut}(X,\mu)$-invariant.
By \cite[Claim~18]{FW} the $\text{Aut}(X,\mu)$-orbit of any free $\Gamma$-action is dense in $\mc A_\Gamma$.
Therefore, it remains to show that $\mc W$ and $\mc S$ contain at least one free action.
Each Bernoullian $\Gamma$-action is free and belongs to $\mc W$.
Since each ergodic action with pure point spectrum has a simple spectrum, $\mc S$ contains a free $\Gamma$-action by Lemma~\ref{MainLem}.
\end{proof}

Let $e_0:=(0,0,1)\in\Gamma$.
Notice that $e_0 (0,j,0) e_0^{-1}=(0,Aj,0)$.
Therefore, $A$ extends naturally to $\Gamma$ via the conjugation by $e_0$.

\begin{Lem}\label{N2N}
Let $\mc H$ be a separable Hilbert space and let $U\colon H \ni h\mapsto U(h) \in \mathcal U(\mc H)$ be a unitary representation of $H$ in $\mc H$.
If $U$ is unitarily equivalent to $U\circ A$ and for each $1\leqslant i<N$, the representation $U\upharpoonright H_i$ of $H_i$ has no non-trivial fixed vector then $\mathcal M(U\upharpoonright F)\subset\{N,2N,\ldots\}\cup\{\infty\}$.
\end{Lem}

\begin{proof}
By the spectral theorem for $U$, there exist a probability measure $\sigma$ on the dual group $\widehat{H}$ and a Borel map $k\colon\widehat{H}\to\N\cup\{\infty\}$ such that the following decomposition holds (up to unitary equivalence):
\begin{equation}\label{Hdecomp}
\mc H = \int_{\widehat{H}}^\oplus \mc H_\omega d\sigma(\omega)
\quad \text{и}\quad
U(h) = \int_{\widehat{H}}^\oplus \omega(h)I_\omega d\sigma(\omega)
\end{equation}
for each $h\in H$, where $\omega\mapsto\mc H_\omega$ is a Borel field of Hilbert spaces, $\dim \mc H_\omega = k(\omega)$ and $I_\omega$ is the identity operator on $\mc H_\omega$.

The inclusion $F\to H$ induces a projection $\pi:\widehat{H}\to\widehat{F}$.
Let $\sigma=\int_{\widehat{F}} \sigma_\chi d\widehat{\sigma}(\chi)$ denote the disintegration of $\sigma$ relative to this projection.
Then we derive from (\ref{Hdecomp}) that
$$
\mc H = \int_{\widehat{F}}^\oplus \mc H'_\chi d\widehat{\sigma}(\chi)
\quad \text{and}\quad
U_g = \int_{\widehat{F}}^\oplus \chi(g)I_\chi d\widehat{\sigma}(\chi),
$$
where $\mc H'_\chi:=\int_{\widehat{H}}^\oplus \mc H_\omega d\sigma_\chi(\omega)$.
Let $l(\chi):=\dim \mc H'_\chi$, $\chi\in\widehat{F}$.
Then
\begin{equation}\label{L}
l(\chi)=\begin{cases}
\sum_{\sigma_\chi(\omega)>0} k(\omega) ,& \text{if $\sigma_\chi$ is purely atomic,}\\
\infty, & \text{otherwise.}
\end{cases}
\end{equation}

Notice that $\sigma\circ A^*$ is the spectral measure for $U\circ A$ where $A^*$ is the dual to $A$ automorphism of $\widehat{H}$, and recall that $A^N=\text{id}$.
Since $U$ is unitarily equivalent to $U\circ A$, it follows that $\sigma$ is equivalent to $\sigma\circ A^*$ and $k=k\circ A^*$.
We may therefore assume that $k$ and $\sigma$ are both invariant under $A^*$ by replacing $\sigma$ with $\frac{1}{N}\sum_{n=0}^{N-1}\sigma\circ {A^*}^n$ if necessary.

We claim that
\begin{equation}\label{OrbLen}
\text{for $\sigma$-almost all $\omega\in\widehat{H}$, $A^*$-orbit of $\omega$ has length $N$. }
\end{equation}
Indeed, otherwise there exist  $1\leqslant m\leqslant N-1$ and $B\subset\widehat{H}$, $\sigma(B)>0$, such that ${A^*}^m\omega = \omega$, $\omega\in B$.
Then by (\ref{Hdecomp}) there exists vector $0\neq\psi\in\mc H$, $\text{supp}\,\psi\subset B$, such that $U(A^m h) \psi = U(h)\psi$ for every $h\in H$.
Equivalently, $U(A^m h-h) \psi = \psi$ for each $h\in H$, i.~e., $\psi$ is a non-trivial fixed vector for $U\upharpoonright H_m$.
However, this contradicts to a condition of the lemma.

Since $Ag=g$ for all $g\in F$, we have $\pi\circ A^*=\pi$.
Therefore, in follows from the invariance of $\sigma$ under $A^*$ that $\sigma_\chi\circ A^*=\sigma_\chi$ for $\widehat{\sigma}$-almost all $\chi\in\widehat{F}$.
Hence (\ref{L}), (\ref{OrbLen}) and the fact $k\circ A^* = k$ imply that  $N\mid l(\chi)$ for $\widehat{\sigma}$-almost all $\chi\in\widehat{F}$, i.~e. $\mathcal M(U\upharpoonright F)\subset\{N,2N,\ldots\}\cup\{\infty\}$.
\end{proof}

Now we are ready to prove the main result of this section (recall that $G$ is isomorphic to $F\subset\Gamma$).

\begin{Th}
For each $T\in \mc S\cap\mc W$, i.~e. for a generic action from $\mc A_\Gamma$, the action $T\upharpoonright F$ is weakly mixing and $\mc M(T\upharpoonright F)=\{N\}$.
\end{Th}

\begin{proof}
Since $T\in\mc W$ and $U_T(e_0) U_T(h) U_T(e_0)^{-1} = U_T(Ah)$ for each $h\in H$, we can apply  Lemma~\ref{N2N} which yields that
$$
\mc M(T\upharpoonright F)\subset\{N,2N,\ldots\}\cup\{\infty\}.
$$
On the other hand, $\mc M(T\upharpoonright\Lambda)=\{1\}$ since $T\in\mc S$.
Recall that $F$ is a subgroup of index $N$ in $\Lambda$.
Hence
$$
\mc M(T\upharpoonright F)\subset\{1,\ldots,N\}.
$$
Therefore, $\mc M(T\upharpoonright F)=\{N\}$.
Since $1\notin \mc M(T\upharpoonright F)$, it follows that
$T\upharpoonright F$ is weakly mixing.
\end{proof}

\section{$\R^p$-actions with homogeneous spectrum} \label{RpSec}

In this section we prove  Theorem~\ref{MainTh} in the case when $G=\Bbb R^p$.

\begin{Th}\label{MainTh2}
Given $p>0$ and $N>1$, there exists a weakly mixing measure preserving $\R^p$-action with homogeneous spectrum of multiplicity $N$.
\end{Th}

From now on we let $J:=\Bbb Z$.
As in the previous section, we denote by $\Gamma$ the semidirect product $\Gamma:=\R^p\times J^\ZZ{N}\rtimes_A(\ZZ{N})$, where $A$ is  the same (as in Section~\ref{DiscSec}) `cyclic' group automorphism.
Consider the following subgroups in $\Gamma$:
\begin{align*}
\Lambda &:=\{ (t,0,i) \mid t\in \R^p, i\in \ZZ{N} \},\\
\Phi &:=\{ (0,j,i) \mid j\in J^\ZZ{N}, i\in \ZZ{N} \},\\
H &:=\{ (t,j,0) \mid t\in \R^p, j\in J^\ZZ{N} \},\\
H_i &:=\{ (0,j,0) \mid j\in J^\ZZ{N}, j(i)=-j(0), j(k)=0 \text{ for } k\neq 0, i  \}\cong J,\quad 1\leqslant i<N,\\
F &:=\Lambda\cap H\cong \R^p.
\end{align*}

We denote by $\mc A_\Gamma$ the set of all measure-preserving $\Gamma$-actions on $X$.
Let $d_0$ stand for a complete metric compatible with the  weak topology on $\text{Aut}(X,\mu)$ \cite{Gl}.
Fix a map $\alpha:\Phi\to (0,1)$ such that $\sum_{f\in \Phi}\alpha(f)<\infty$ and equip $\mc A_\Gamma$ with the metric
$$
d(T,S):=\sup_{t\in [-1,1]^p}d_0(T_{(t,0,0)},S_{(t,0,0)})+\sum_{f\in \Phi}\alpha_fd_0(T_f,S_f).
$$
It is easy to see that $d$ is compatible with the topology of uniform convergence on the compact subsets in $\Gamma$.
Since $(\text{Aut}(X,\mu),d_0)$ is a Polish group, one can deduce easily that $(\mc A_\Gamma,d)$ is a Polish space.
Moreover, the natural action of Aut$(X,\mu)$  on $\mc A_\Gamma$ by conjugation is continuous.
In a similar way, given any locally compact second countable group $\Lambda$, a Polish topology is introduced on the set $\mc A_\Lambda$ of
$\mu$-reserving actions of $\Lambda$.

In the previous section we widely used the following facts about the space of actions of a discrete countable group, say $\Sigma$:
\begin{enumerate}
\item[(a)] the set of free $\Sigma$-actions is $G_\delta$ in $\mc A_\Sigma$ \cite{GlK},
\item[(b)] the Aut$(X,\mu)$-orbit of any free $\Sigma$-action is dense in $\mc A_\Sigma$ \cite{FW}.
\end{enumerate}

We first briefly demonstrate similar  facts for actions of some
{\it continuous} locally compact groups (because we were unable to  find them in the literature).

\begin{Lem}\label{lab}
Given a locally compact Polish group $\Lambda$, the subset $\mc F$ of free $\Lambda$-actions  is $G_\delta$ in $\mc A_\Lambda$.
\end{Lem}

\begin{proof}
Fix a sequence  $\mc A_i=(A_1^{(i)},\dots,A_{2^i}^{(i)})$ of finite partitions of $X$ such that $\mc A_{i+1}$ refines $\mc A_i$ and $\bigvee_{i>0}\mc A_i=\fr B$.
Select a sequence of compacts $K_1\subset K_2\subset\cdots$ in $\Lambda$ such that $\bigcup_{i>0}K_i=\Lambda\setminus\{1_\lambda\}$.
Then it remains to note that
$$
\mc F=\bigcap_{l=1}^\infty\bigcap_{m=1}^\infty\bigcap_{n=1}^\infty \bigcup_{i=0}^\infty\{T\in\mc A_\Lambda\mid \sum_{j\in\mc J_{i,l,n}}\mu(A^{(i)}_j)>1-\frac 1m\},
$$
where $\mc J_{i,l,n}=\{j\mid \max_{k\in K_l}\mu(T_kA^{(i)}_j\cap A^{(i)}_j)<\frac 1n\mu(A^{(i)}_j)\}$.
\end{proof}

Now we return to the group $\Gamma$ and note that it is monotilable.
This means that there exist a F{\o}lner sequence $\{F_n\}_{n\in\N}$ and a sequence $\{C_n\}_{n\in\N}$ of countable subsets in $\Gamma$ such that $\{F_n c\mid c\in C_n\}$ is a partition of $\Gamma$ for each $n$.
(We will also assume that $F_1\subset F_2\subset\cdots$ and $\bigcup_{n>0}F_n=\Gamma$.)
 Hence a ``classical'' Rokhlin lemma holds for ergodic free actions of $\Gamma$ as follows.

\begin{Lem}[\cite{Se}, \cite{OW}]\label{Rokhlin}
Let $T$ be an ergodic free action of $\Gamma$.
Given $\varepsilon > 0$ and
$n>0$, there exist a subset $A\subset X$ with $\mu(A)>1-\varepsilon$, a standard probability space $(Y,\nu)$ and a one-to-one measure preserving mapping $\varphi\colon (A,\mu) \to (Y\times F_n,\nu\times\lambda)$ such that
$$
\varphi T_g \varphi^{-1}(y,h) = (y,gh)
\quad\text{for}\quad  (y,g)\in Y\times F_n, h\in F_n\cap(g^{-1}F_n),
$$
where $\lambda$ is the  Haar measure on $\Gamma$ normed in such a way that $\lambda(F_n)=\mu(A)$.
\end{Lem}

We will call $(A,\varphi,Y)$ an \emph{$(F_n,\varepsilon)$-Rokhlin tower}.
The following assertion is a natural corollary of Lemma~\ref{Rokhlin}.

\begin{Lem}\label{RpFW}
The $\text{Aut}(X,\mu)$-orbit of any $\Gamma$-action $T\in\mc E$ is dense in $\mc E$.
\end{Lem}

\begin{proof}
Let $T,T'\in\mc E$.
Fix  a Haar measure $\lambda$ on $\Gamma$.
Take $\varepsilon > 0$, a compact subset $K\subset \Gamma$ and a finite sequence $B_1,\ldots, B_k$ of measurable subsets in $X$.
Take $n$ large so that $K\subset F_n$ and $\lambda(F_n\bigtriangleup(K^{-1}F_n))<\varepsilon\lambda(F_n)$.
By  Lemma \ref{Rokhlin}, there exist $(F_n,\varepsilon)$-Rokhlin towers $(A,\varphi,Y)$ and $(A',\varphi',Y')$ for $T$ and $T'$, respectively.
Without loss of generality we may assume that $\mu(A) = \mu(A')$.
Let $S\colon Y \to Y'$ be a measure preserving bijection.
Then we define $R\colon A\to A'$ by setting $R := (\varphi')^{-1}(S\times I)\varphi$.
It follows that $R T_g x = T'_g R x$ for all $g\in F_n$, $x\in \varphi^{-1}(Y\times(F_n\cap g^{-1}F_n))\subset A$.
Extend $R$ to a one-to-one measure preserving transformation of $X$ in an arbitrary way.
The inclusion $RT_gR^{-1} B_i \bigtriangleup T'_g B_i \subset \{ x\in A' \mid R T_g x \neq T'_g R x \}$
implies that
$$
\mu(R T_g R^{-1} B_i \bigtriangleup T'_g B_i)<2\varepsilon
$$
for all $g\in K$, $1\leqslant i\leqslant k$. This implies that $T'$ belongs to the closure of the $\text{Aut}(X,\mu)$-orbit of $T$.
\end{proof}

Thus we established analogues of (a) and (b) for continuous groups.
Now let $\mc E \subset \mc A_\Gamma$ stand for the subset of all free actions $T$ such that
$T\upharpoonright\Lambda$ is ergodic.

\begin{Lem}\label{EisGd}
$\mc E$ is $G_\delta$  in $\mc A_\Gamma$.
\end{Lem}

\begin{proof}
Since the map $\mc A_\Gamma\ni T\mapsto T\upharpoonright \Lambda \in\mc A_{\Lambda}$ is  continuous, it suffices to show that
the set of  ergodic $\Lambda$-actions is $G_\delta$ in $\mc A_\Lambda$ (and apply also Lemma~\ref{lab}).
Since the von Neumann mean ergodic theorem holds for the $\Lambda$-actions along a fixed F{\o}lner sequence, we deduce that the set of ergodic actions is  $G_\delta$ in $\mc A_\Lambda$.
For that, we argue almost verbally as in the proof of Lemma~\ref{WM}.
\end{proof}

It follows that $\mc E$ is  Polish  when endowed with the induced topology.

We now state an analogue of Lemma~\ref{MainLem}.

\begin{Lem}\label{RpErgPPS}
There exists an action $T\in\mc E$
such that the restriction $T\restriction\Lambda$ has a pure point spectrum.
\end{Lem}

\begin{proof}

Take rationally independent  irrationals $\alpha_{n}$, $n=1,\ldots,p+N$.
Define one-to-one homomorphisms $R\colon \R^p\to\T^{p+1}$ and $S_k\colon J\to\T^{p+1}$, $0\leqslant k<N$, by setting
\begin{align*}
R(t) &:=(e^{2\pi i t_1},\ldots,e^{2\pi i t_p},e^{2\pi i(\alpha_{1} t_1+\cdots+\alpha_{p} t_p)})\\
S_k(j) &:= (1,\dots,1,e^{2\pi i\alpha_{p+k+1}j}),
\end{align*}
$t=(t_1,\ldots,t_p)\in\R^p$, $j\in J$.
Then the subgroups $R(\R^p), S_0(J),\ldots,S_{N-1}(J)$ are independent in $\T^{p+1}$.
We set $X:=\T^{p+1}\times (\ZZ{N})$.
Then $X$ is a compact Abelian group.
We equip $X$ with Haar measure.
Now we define an action $T$ on $X$ as follows:
$$
T_{(t,j,m)}(z,n):=(z+R(t)+S_{n+m}(j(n+m)),n+m),
$$
$(t,j,m)\in\Gamma$, $(z,n)\in X$.
Obviously, $T$ is free.
Since $T_{(t,0,n)}x=x+(R(t),n)$, $(t,0,n)\in\Lambda$, $x\in X$, and the subgroup $\{ (R(t),n) \mid (t,0,n)\in\Lambda \}$ is dense in $X$,
 it follows that the action $T\upharpoonright\Lambda$ is ergodic with a pure point spectrum.
\end{proof}

The following assertion is a direct analogue of Lemma~\ref{WSresidual1}.

\begin{Lem}
The following two subsets are residual in $\mc E$:
\begin{align*}
\mc W &:=\{ T\in\mc E \mid T\upharpoonright H_i \text{ is weakly mixing for each $1\leqslant i < N$} \}
\quad\text{and}\\
\mc S &:= \{ T\in\mc E \mid T\upharpoonright\Lambda \text{ has a simple spectrum} \}
\end{align*}
\end{Lem}

\begin{proof}
The subset of weakly mixing $H_i$-actions and the subset of $\Lambda$-actions with simple spectrum are $G_\delta$ in $\mc A_{H_i}$ and $\mc A_\Lambda$ respectively (the proof is  similar to the proofs of Lemmata~\ref{WM}, \ref{SimSp}).
Therefore, $\mc W$ and $\mc S$ are both $G_\delta$ in $\mc E$.
Notice that $\mc W$ and $\mc S$ are $\text{Aut}(X,\mu)$-invariant.
Therefore in view of Lemma~\ref{RpFW} it remains to show that $\mc W$ and $\mc S$ are non-empty.
 Consider an action of $\Gamma$ on itself  by translations.
This action preserves the ($\sigma$-finite, infinite) Haar measure.
The corresponding Poisson suspension of this action is  a probability preserving free $\Gamma$-action and it belongs to $\mc W$ (see \cite{OW}).
It is an analogue of Bernoulli action for continuous groups.
On the other hand, $\mc S\ne\emptyset$ by Lemma~\ref{RpErgPPS}.
\end{proof}

Lemma \ref{N2N} extends to the setting of continuous groups almost verbally.
Therefore  we conclude the following

\begin{Th}
For each $T\in \mc S\cap\mc W$, i.~e. for a generic action from $\mc E$, the action $T\upharpoonright F$ is weakly mixing and $\mc M(T\upharpoonright F)=\{N\}$.
\end{Th}

Recall that $F$ is isomorphic to $\R^p$, and hence  Theorem~\ref{MainTh2} is shown.

\section{General case}

Our purpose in this section is to prove Theorem~\ref{MainTh}
 in the full generality.
 Thus $G=\Bbb R^p\times G'$ for some $p\ge 0$, where $G'$ is a locally compact Polish group which contains an open compact subgroup $G'_0$.
As in the two previous sections, we select a certain  discrete countable Abelian subgroup $J$ (depending on $G$) and consider a semidirect product
$\Gamma=G\times J^{\Bbb Z/N\Bbb Z}\rtimes_A\Bbb Z/N\Bbb Z$.
The set $\mc A_\Gamma$ of probability preserving $\Gamma$-actions is equipped with the Polish topology of uniform convergence on the compact subsets of $\Gamma$.
Arguing as in Sections~\ref{DiscSec} and \ref{RpSec} one can show that for a generic action $T$ from $\mc A_\Gamma$, the restriction of $T$ to $G$   has a homogeneous spectrum of multiplicity $N$.
It remains only to explain how to select $J$ and how to construct a free $\Gamma$-action whose restriction to the subgroup $G\times\{0\}\times \Bbb Z/N\Bbb Z$ is ergodic with a pure point spectrum.

(i) Let $p>0$.
Then we set $J=\Bbb Z$.
Let $T$ be a free action of $\Bbb R^p\times J^{\Bbb Z/N\Bbb Z}\rtimes_A\Bbb Z/N\Bbb Z$ such that $T\restriction(\Bbb R^p\times \{0\}\times\Bbb Z/N\Bbb Z)$ is ergodic with a pure point spectrum.
Such an action is constructed in Lemma~\ref{RpErgPPS}.
Now take any free ergodic $G'$-action $T'$ with a pure point spectrum.
Then the product action  $T\otimes T'$ (considered as a $\Gamma$-action)
is as desired.

(ii) Let $p=0$ but $G'=G'_0\times G'/G'_0$.
Since $G$ is non-compact,  the quotient group $G'/G'_0$ is infinite.
Then we select $J$ and construct  a free action $T$ of the corresponding group $G'/G'_0\times J^{\Bbb Z/N\Bbb Z}\rtimes_A\Bbb Z/N\Bbb Z$ such that $T\restriction(G'/G'_0\times \{0\}\times\Bbb Z/N\Bbb Z)$ is ergodic with a pure point spectrum as in Lemma~\ref{MainLem}.
Let $T'$ be an action of $G_0'$ on itself (equipped with Haar measure) by translations.
Then the product action $T'\times T$ is as desired.

(iii) Suppose that $p=0$, $G_0'$ is not a direct summand in $G'$ and
 there is no $k>0$ such that  $k\cdot g=0$ for all $g\in G'/G_0'$.
Let $K$ be a compact Abelian group and $R:G'\to K$ be a continuous one-to-one group homomorphism with $\overline{R(G')}=K$.
Of course, $R(G_0')$ is closed in $K$.
Moreover, a one-to-one homomorphism
$$
R': G'/G'_0\ni g'+G'_0\mapsto R(g')+R(G'_0)\in K/R(G'_0)
$$
is well defined.
It is easy to see that  $R'(G'/G'_0)$ is dense in $K/R(G'_0)$.
It follows that $K/R(G'_0)$ has elements of infinite order (\cite[Theorem 25.9]{HR}).
Hence arguing as in the proof of Lemma~\ref{MainLem}, we set $J=\Bbb Z$ and define group homomorphisms
$S_0',\dots, S_{N-1}':J\to K/R(G'_0)$ in such a way that
the subgroups $R'(G'/G'_0)$, $S_0'(J),\dots, S_{N-1}'(J)$ are independent in
$K/R(G'_0)$.
Then there are homomorphisms $S_0,\dots, S_{N-1}:J\to K$ such that
$S_i(k)+R(G'_0)=S_i'(k)$, $k\in J$.
The subgroups
$R'(G')$, $S_0(J),\dots, S_{N-1}(J)$ are independent in $K$.
Hence we can define the sought-for action of $\Gamma$ on $K\times \Bbb Z/N\Bbb Z$ in a similar way as in the proof of Lemma~\ref{MainLem}.

The problem why this does not work if the orders of all elements of $G'/G'_0$ are bounded from above is that then $J$ is a torsion group and we, in general, can not `cover' a homomorphism $S':J\to K/R(G'_0)$ with
 a homomorphism $S:J\to K$.



\end{document}